\documentclass[11pt]{article}

\usepackage{amsmath}
\usepackage{amssymb}
\usepackage[PostScript=dvips]{diagrams}

\DeclareMathAlphabet{\scr}{U}{eus}{m}{n}

\newcommand\N{{\mathbb N}}
\newcommand\Z{{\mathbb Z}}
\newcommand\Q{{\mathbb Q}}
\newcommand\F{{\mathbb F}}
\newcommand\PP{{\mathbb P}}
\newcommand\A{{\mathbb A}}
\newcommand\G{{\mathbb G}}

\newcommand\fp{{\mathfrak p}}
\newcommand\fm{{\mathfrak m}}

\newcommand\cA{{\cal A}}
\newcommand\cB{{\cal B}}
\newcommand\cC{{\cal C}}
\newcommand\cF{{\cal F}}

\newcommand\cN{{\cal N}}
\newcommand\cO{{\cal O}}
\newcommand\cP{{\cal P}}
\newcommand\cR{{\cal R}}
\newcommand\cT{{\cal T}}

\newcommand\cNN{\widehat{\cal N}^-}

\newcommand\sL{{\scr L}} 
\newcommand\sM{{\scr M}}

\newcommand\Cnorm{\widetilde{C}}

\newcommand\Mloc{{\rm\bf M}^{\rm loc}}
\newcommand\Mlocs{\overline{\rm\bf M}^{\rm loc}} 
\newcommand\Mloct{\widetilde{\rm \bf M}^{\rm loc}}

\newcommand\xt{[\hspace{-.12em}[ t ]\hspace{-.12em} ] } 
\newcommand\xT{(\!( t )\!)}

\newcommand\tens\otimes

\newcommand\lto{\longrightarrow}

\newcommand{\Spec}{\mathop{\rm Spec}}
\newcommand{\diag}{\mathop{\rm diag}}
\newcommand{\Grass}{\mathop{\rm Grass}\nolimits}

\newcommand\qed{\hfill$\square$}

\newtheorem{thm}{Theorem}[section]
\newtheorem{stz}[thm]{Proposition}
\newtheorem{lem}[thm]{Lemma}
\newtheorem{Def}[thm]{Definition}
\newtheorem{kor}[thm]{Corollary}

\begin{document}

\title{On the flatness of local models\\ for the symplectic group} 
\author{Ulrich G\"ortz}
\date{October 2000}
\maketitle

\section{Introduction}

In order to study arithmetic properties of a variety over
an algebraic number field, it is desirable to have a model
over the ring of integers. We are interested in the
particular case of a Shimura
variety (of PEL type), and ask for a model over $O_E$,
where $E$ is the completion of the reflex field at
a prime of residue characteristic $p>0$.
Since a Shimura variety of PEL type is a moduli
space of abelian varieties with certain additional structure, 
it is natural to define
a model by posing the moduli problem over $O_E$; in the case of
a parahoric level structure at $p$, which is the case we are
interested in, such models have been defined by Rapoport 
and Zink \cite{RZ}.

These models are almost never smooth, and it is interesting to study
the singularities occuring in the special fibre. It even 
is not obvious whether the models are flat over $O_E$, 
although flatness certainly belongs to the minimum requirements
for a reasonable model.
Rapoport and Zink have conjectured that their models
are indeed flat.
In the case of a unitary group that splits over an 
unramified extension of $\Q_p$, flatness has been proved
in \cite{G}. In this article we will extend this result
to the case of the symplectic group, i.e. we will show the flatness 
for groups of the form ${\rm Res}_{F/\Q_p} GSp_{2r}$,
where $F$ is a finite {\em unramified} extension of $\Q_p$.
In the case where the underlying group does not
split over an unramified extension, the flatness conjecture
is not true as it stands, as has been pointed out by Pappas \cite{P};
see also \cite{PR} for ideas on how to proceed in this case.

Several special cases of our main theorem have been known before,
by the work of Chai and Norman
\cite{CN}, de Jong \cite{dJ} and Deligne and Pappas \cite{DP}.

To study local questions such as flatness, it is useful
to work with the so called local model which \'etale locally
around each point of the special fibre coincides with the model
of the Shimura variety, but which can be defined in terms
of linear algebra and is thus easier to handle (cf. \cite{dJ},  
\cite{RZ}).

The proof of the flatness conjecture in the symplectic case then
starts with the same approach as in the unitary case. Namely, we embed
the special fibre of the local model into the affine flag variety, and
identify it with the intersection of certain Schubert varieties,
coming from simpler local models. Via the theory of Frobenius
splittings we can conclude that this intersection is reduced.
To make the proof work some additional ingredients are needed,
most notably Faltings' result on the normality of Schubert varieties
in the affine flag variety (see \cite{F} and section \ref{normality} 
below). Furthermore, to analyse the 'simplest'
singular local models, we use a theorem of de Concini
about so-called doubly symplectic standard tableaux (see \cite{dC}
and section \ref{deConcini}). 
More precisely, an open neighborhood of their worst singularity 
in the special fibre is essentially isomorphic to $\Spec R$, where
$$  R= \F_p[c_{\mu \nu}, \mu, \nu = 1, \dots 2i] / 
    (CJ_{2i}C^t, C^tJ_{2i}C).
$$
Here $C = (c_{\mu\nu})_{\mu,\nu}$ and  
 $J_{2i} =\left( \begin{array}{cc} & {\mathbf J} \\ -{\mathbf J} & 
\end{array} \right)$, where 
 ${\mathbf J}$ is the $i \times i$-matrix with
1's on the second diagonal and 0's elsewhere. 
In loc. cit., de Concini gives a $\F_p$-basis of $R$, and this will
enable us to show that $R$ is reduced (see proposition \ref{Rred}).
Note that the same singularity has also been studied
by Faltings \cite{F1}.

The case treated by Chai and Norman respectively by Deligne and Pappas
corresponds to the equations
$$ CD = DC = p,\quad C = C^t,\quad D = D^t, $$
where $C$ and $D$ are $r \times r$-matrices of indeterminates
over $\Z_p$. They show that these singularities are Cohen-Macaulay,
and that the special fibre of the corresponding scheme is
reduced. The proofs rely on the theory of algebras
with straightening laws. Our approach gives a new proof for this
result and it may be interesting to note that for this special
case we do not need to invoke de Concini's result. See remark i)
after theorem \ref{mainthm}.

The precise statement of the main theorem as well as a sketch of
its proof can be found in section \ref{statement}.
The
proof proper of the main theorem is actually rather short. However, to
make the paper more accessible, we have included an exposition of some
of the results of others that we use. In section \ref{normality}
 we give an exposition
of Faltings'  proof of the normality result mentioned above in the case
of interest to us, and in section \ref{deConcini}
 we explain the result of de Concini
and the crucial consequence of it that we need. The final section 
\ref{proof}
brings everything together and contains the proof of the main results.

Finally it is a pleasure to express my gratitude to those who have helped
me with this work. First of all, I am much obliged to M. Rapoport for
many helpful discussions and for his steady interest in my work.
I am grateful to S. Orlik and T. Wedhorn
for several useful conversations; T. Wedhorn also made several useful 
remarks on this text.
Furthermore I would like to thank G. Faltings for briefly explaining to me
his proof of the theorem on the normality of affine Schubert 
varieties.
This manuscript has been finished during my stay at the Institute for Advanced
Study. I am grateful to the Institute for providing a great 
working environment and to the Deutsche Forschungsgemeinschaft
and the National Science Foundation (grant DMS 97-29992) for their support.

\section{The standard local model for the symplectic group}
\label{statement}

First of all, let us give the definition of the
 {\em standard local model} which is the object of
our studies.
Let $O$ be a complete discrete valuation ring with perfect 
residue class field of characteristic $p > 0$.
Let $\pi$ be a uniformizer of $O$ and let $k$ be an algebraic
closure of the residue class field of $O$.

Fix a positive integer $r$, and let $n = 2r$.

Denote the quotient field of $O$ by $K$. Let $e_1, \dots, e_n$
be the canonical basis of $K^n$. We endow $K^n$ with the standard
symplectic form, 
given by the matrix 
 $J_{2r} =\left( \begin{array}{cc} & {\mathbf J} \\ -{\mathbf J} & 
\end{array} \right)$, where 
 ${\mathbf J} = (\delta_{r-i+1, j})_{ij}$ is the $r \times r$-matrix with
1's on the second diagonal and 0's elsewhere. 

Let $\Lambda_i$, $0 \le i \le n-1$, be the free $O$-module of rank 
$n$ with basis 
 $\pi^{-1} e_1, \dots, \pi^{-1} e_{i}, e_{i+1}, 
   \dots, e_n$.
This gives rise to a complete lattice chain
\begin{diagram}
\cdots & \rTo & \Lambda_0 & \rTo & \Lambda_1 
& \rTo & \cdots & \rTo & \Lambda_{n-1} & \rTo & \pi^{-1}\Lambda_0 & \cdots
\end{diagram}

Furthermore, the lattice chain is selfdual, i.e. for each $\Lambda$
occuring in the lattice chain, the dual
$$ \Lambda^\ast = \{ x \in K^n; \langle x, y \rangle \in O \text{ for all }
                                   y \in \Lambda \} $$
appears as well. More precisely, $\Lambda_0^\ast = \Lambda_0$,
 $\Lambda_r^\ast = \pi \Lambda_r$, and 
 $\Lambda_i^\ast = \pi \Lambda_{2r-1}$, $i = 1, \dots, r-1$.

Finally, choose a subset $I = \{ i_0 < \cdots < i_{m-1} \}$,
such that for each $i \in I$, $1 \le i \le 2r-1$, also $2r-i \in I$.

The standard local model $\Mloc_I$, as defined by Rapoport
and Zink, is the $O$-scheme 
representing the following functor (cf. \cite{RZ}, definition 3.27):

For every $O$-scheme $S$, $\Mloc_I(S)$ is the set of isomorphism classes
of commutative diagrams 
\begin{diagram}
\Lambda_{i_0, S} & \rTo & \Lambda_{i_1,S} & \rTo & \cdots 
& \rTo & \Lambda_{i_{m-1},S} & \rTo^\pi & \Lambda_{i_0,S} \\
\uInto              &      & \uInto        &      &   &&\uInto & &\uInto   \\
\cF_0     & \rTo & \cF_1  & \rTo & \cdots & \rTo & \cF_{m-1}
  &\rTo & \cF_0
\end{diagram}
where $\Lambda_{i,S}$ is $\Lambda_i \tens_{O} \cO_S$, and where 
the $\cF_\kappa$ are locally free $\cO_S$-submodules of rank $r$ which
Zariski-locally on $S$ are direct summands of $\Lambda_{i_\kappa,S}$. 

Furthermore, the $\cF_i$ have to satisfy the following duality condition:
for each $i \in I$, the map
$$ \cF_i \lto \Lambda_{i,R} \cong \widehat{\Lambda}_{2r-i,R} \lto 
   \widehat{\cF}_{2r-i} $$
is the zero map.

It is clear that this functor is indeed representable. In fact, 
 $\Mloc_I$ is a
closed subscheme of a product of Grassmannians.
More precisely, the local model for the symplectic group $GSp_{2n}$ is
a closed subscheme of the local model associated to $GL_{2n}$ (for the same
 subset $I$); for its definition see \cite{RZ} or \cite{G}.

We write $\Mloc := \Mloc_{\{0, \dots, 2r-1\}}$.
Furthermore, we will often write $\Mloc_{i_0, \dots, i_{m-1}}$
instead of $\Mloc_{\{i_0, \dots, i_{m-1} \}}$.

The aim of this article is to prove the following theorem:

\begin{thm} \label{mainthm}
The local model $\Mloc_I$  is flat over $O$, 
and its special fibre is reduced.
The irreducible components of the special fibre
are normal with rational singularities,
so in particular they are Cohen-Macaulay.
\end{thm}

Sketch of proof (in the Iwahori case $I = \{ 0, \dots, 2r-1 \}$):

We embed the special fibre $\Mlocs$ of the local model 
into the affine flag variety
$\cF = Sp_{2r}(k\xT)/B$. Furthermore, we denote by 
 $\Mloct_I$ the inverse image of $\Mlocs_I \subset Sp_{2r}(k\xT)/P^I$
under the projection $\cF \lto Sp_{2r}(k\xT)/P^I$.

Then we have
$$
 \Mlocs = \Mloct_0 \cap \Mloct_r \cap 
           \bigcap_{i=1}^{r-1} \Mloct_{i, 2r-i}.
$$

Now the proof of the flatness 
theorem can be divided into the following three steps.

\begin{itemize}

\item The special fibres of the local models $\Mloc_0$, $\Mloc_r$, 
      $\Mloc_{i, 2r-i}$, $ i= 1, \dots, r-1$,
       are reduced (and irreducible) and thus are Schubert varieties
      (see section \ref{simplest}).
      Here we use a result of de Concini, which is explained in
      section \ref{deConcini}.

\item Schubert varieties in the affine flag variety 
      for the symplectic group are normal
      and Frobenius split, as was shown by Faltings (see section
      \ref{normality}). This implies 
      that intersections
      of Schubert varieties are reduced. In particular we see
      that $\Mlocs$ is reduced (see section \ref{inters}). 

\item The generic points of the irreducible components of the special
fibre $\Mlocs$ can be lifted to the generic fibre (see section
\ref{liftability}). Since $\Mlocs$
is reduced, this implies that $\Mloc$ is flat over $O$.

\end{itemize}

{\bf Remarks.}
i) Choosing $I = \{ 0, r \}$, we get the result of Chai
and Norman \cite{CN} respectively Deligne
and Pappas \cite{DP}. Since in this case 
we have 
$$
 \Mloct_I = \Mloct_0 \cap \Mloct_r,
$$
and $\Mloct_0$ and $\Mloct_r$ are smooth, we can conclude
that $\Mloct_I$ (and thus $\Mlocs_I$) is reduced without
using de Concini's result. Thus in this case our approach gives a 
proof that does not use the theory of algebras 
with straightening laws, unlike those of Chai/Norman and
Deligne/Pappas.

ii) While the theorem as it stands deals with the case of $GSp_{2r}$
over $\Q_p$, it can easily be generalized to the case of groups
of the form ${\rm Res}_{F/\Q_p} GSp_{2r}$, where $F$ is an unramified 
extension of $\Q_p$. See \cite{RZ}, chapter 3, for the definition
of the local model in this case; see \cite{G}, section 4.6, for details 
on how to obtain the generalization.

iii) On the other hand, by some kind of continuity argument, it can
be shown that the theorem holds as well if the characteristic
of the residue class field of $O$ is 0. Confer \cite{G}, section 4.5.

iv) It is an interesting question, if $\Mloc$ itself is Cohen-Macaulay.
In view of the flatness, this is equivalent to $\Mlocs$ being
Cohen-Macaulay. This question can be reduced to a purely combinatorial
problem in the affine Weyl group; unfortunately the combinatorics
involved seems to be rather difficult.

\section{The affine flag variety for the symplectic group}
\label{normality}

\subsection{Definition}

In this section we will give the definition of
the affine flag variety and reproduce Faltings' proof for the normality
of Schubert varieties contained in it. 
We will mostly follow Faltings' paper \cite{F}, and no
claim for originality is made. 
Furthermore, we restrict ourselves to the
case of the symplectic group.

Let $G= Sp_{2r}$ (over $\Z$). 
All of the following obviously works for $SL_n$ as well,
and in fact even works for any split, simply connected, semi-simple and
simple group.
We denote by $T$ the torus of diagonal matrices in $G$,
and by $B$ the Borel subgroup of upper triangular matrices.

The loop group $LG$ is the ind-scheme (over $\Z$) such that
for all rings $R$
 $$ LG(R) = G(R\xT). $$
Furthermore, we have the (infinite dimensional) schemes $L^+G$
and $L^{++}G$,
\begin{eqnarray*}
L^+G(R) & = & G(R\xt) \\
L^{++}G & = & \text{kernel of the reduction map } L^+G \lto G,
\end{eqnarray*}
and the ind-schemes
 $ L^-G$ and $L^{--}G$, given by
\begin{eqnarray*}
L^-G(R) & = & G(R[t^{-1}]) \\
L^{--}G & = & \text{kernel of the reduction map } L^-G \lto G.
\end{eqnarray*}
Inside $L^+G$, we have the standard Iwahori subgroup $\cB$,
which is the
inverse image of the standard Borel subgroup of $G$
under the reduction map $L^+G \lto G$.

The affine flag variety $\cF$ is, by definition,
 the quotient $LG / \cB$ in the sense of $fpqc$-sheaves.

Just as the usual flag variety can be interpreted as the
space of flags in a vector space, the affine flag variety
can be seen as a space of lattice chains.
Let us make this more explicit.

\begin{Def} Let $R$ be a ring. A lattice in $R\xT^n$
is a sub-$R\xt$-module $\sL$ of $R\xT^n$ which is locally free of
rank $n$, and such that $\sL \tens_{R\xt} R\xT = R\xT^n$.

A sequence $\sL_0 \subset \sL_1 \subset \dots \subset
  \sL_{n-1} \subset t^{-1} \sL_0$
of lattices in $R\xT^n$ is called a complete lattice chain, if 
 $\sL_{i+1}/\sL_i$ is a locally free $R$-module of rank 1 for all
 $i$.
\end{Def}

By adding all $t^N \sL_i$, $N \in \Z$, to a lattice chain $(\sL_i)_i$,
we get a complete {\em periodic} lattice chain. 

We endow $R\xT^n$ with the alternating pairing 
$\langle \cdot , \cdot \rangle$ given by the matrix
$J_{2r} = \left( \begin{array}{cc} & {\mathbf J} \\ -{\mathbf J} & 
\end{array} \right)$, where 
 ${\mathbf J} = (\delta_{r-i+1, j})_{ij}$ is the $r \times r$-matrix with
1's on the second diagonal and 0's elsewhere.

Then we call a lattice chain selfdual, if for all
lattices $\sL$ in the lattice chain the dual
$$\sL^\ast = \{ x \in K^n; \langle x, y \rangle \in R\xt \text{ for all }
                                   y \in \sL \}
$$
also occurs in the periodic lattice chain.

For example, the standard lattice chain
$$ \lambda_i = t^{-1} R\xt^i \oplus R\xt^{n-i},\quad i=0, \dots, n-1, $$
is selfdual; more precisely we have  $\lambda_0^\ast = \lambda_0$,
 $\lambda_r^\ast = t \lambda_r$, and 
 $\lambda_i^\ast = t \lambda_{2r-1}$, $i = 1, \dots, r-1$.

Finally, we call a lattice chain $(\sL_i)_i$ special, if
 $\bigwedge^n \sL_0 = R\xt$ inside $\bigwedge^n R\xT^n = R\xT$.
For a selfdual lattice chain, this is equivalent to the condition
that $\sL_0$ is a selfdual lattice.

\begin{stz} \label{F_latticechains}
We have a functorial isomorphism
\begin{diagram}
\cF(R) & \rTo^\cong & \{ (\sL_i)_i \text{ special selfdual 
                                complete lattice chain in } R\xT^n \} 
\end{diagram}
In particular, $\cF$ is an ind-scheme again.
\end{stz}

{\em Proof.} The map is of course the map induced by
$$ g \mapsto g \cdot (\lambda_i)_i. $$
It can be seen that the resulting morphism is surjective by
proving that special selfdual complete lattice chains
locally on $R$
admit a normal form. See the appendix to chapter 3 in \cite{RZ}.
\qed

We get a variant of the proposition by considering not special
lattice chains, but $r$-special lattice chains, i. e. lattice
chains $(\sL_i)_i$ such that $\bigwedge^n \sL_0 = t^r R\xt$.
This will be useful for embedding the special fibre of the local 
model into $\cF$.

The affine Weyl group $W_a$ for the symplectic group
 is the semidirect product 
$$ W_r \ltimes \{ (x_1, \dots, x_n, -x_n, \dots, -x_1);
                       x_i \in \Z \}, $$
where $W_r = {\mathfrak S}_r \ltimes \{ \pm 1  \}^r$
is the finite Weyl group of the symplectic group.
The affine Weyl group $W_a$ for $Sp_{2r}$ is a subgroup
of the affine Weyl group for $SL_{2r}$ in a natural way.

Denote the simple reflections in $W$ by $s_1, \dots, s_r$;
these together with the affine simple reflection $s_0$
generate $W_a$. 

Furthermore, let $\cP_i$ be the parahoric subgroup of $LG$ 
generated by $\cB$ and $s_i$. More generally, for $I \subseteq
 \{ 0,\dots, r \}$, we write $\cP_I$ for the parahoric subgroup
generated by $\cB$ and $s_i$, $i \in I$, and 
 $\cP^I = \cP_{  \{ 0,\dots, r \} - I }$.

Corresponding to $s_i$, $i = 0, \dots, r$, we have a subgroup
 $SL_2$ inside $\cP_i$.
Concretely, the subgroup $SL_2(R)$ in $LG(R)$ corresponding to $s_i$, 
$i\in \{1, \dots, r-1\}$ consists of the matrices
$$ \diag(1^{i-1}, \left( \begin{array}{cc} a & b \\ c & d 
   \end{array}\right), 1^{2r-2i}, \left( \begin{array}{cc} a & -b \\ -c & d 
   \end{array}\right),  1^{i-1}), \quad 
\left( \begin{array}{cc} a & b \\ c & d 
   \end{array}\right) \in SL_2(R), 
$$
and the subgroup corresponding to $s_r$ is 
$$ \diag(1^{r-1}, \left( \begin{array}{cc} a & b \\ c & d 
   \end{array}\right), 1^{r-1}), \quad 
\left( \begin{array}{cc} a & b \\ c & d 
   \end{array}\right) \in SL_2(R).
$$
The subgroup $SL_2$ corresponding to the affine simple reflection
 $s_0$ is
$$ \left( \begin{array}{ccccc} 
a  &   &        &   & t^{-1}b \\
   & 1 &        &   &         \\
   &   & \ddots &   &         \\
   &   &        & 1 &         \\
tc &   &        &   &  d      
\end{array} \right), \quad 
\left( \begin{array}{cc} a & b \\ c & d 
   \end{array}\right) \in SL_2(R).
$$

The unipotent subgroup $\left\{ \left( \begin{array}{cc} 1 & b \\ 0 & 1 
          \end{array} \right) \right\} $ of $SL_2$
gives us, for each affine simple root $\alpha$, a subgroup 
 $\G_a \cong U_\alpha \subset LG$, and by conjugation with elements
of the affine Weyl group we get a subgroup $U_\alpha$ for
every affine root $\alpha$.

Now let $k$ be an algebraically closed field.
We have the Bruhat decomposition
$$ \cF (k) = \bigcup_{w \in W_{a}} \cB w \cB / \cB.$$
The closure $C_w$ of $\cB w \cB/\cB$ is called an
(affine) Schubert variety. It is a projective algebraic variety.
We have $C_{w'} \subseteq C_w$ if and only if $w' \le w$
with respect to the Bruhat order.

In the case of characteristic 0, the normality of Schubert varieties
has been known for some time.  Indeed Faltings has proved (see \cite{BL},
\cite{LS})
that in characteristic 0 the affine flag variety defined above
coincides with the one arising in the theory
of Kac-Moody algebras. 
This 'Kac-Moody affine flag variety' has been investigated by
Kumar \cite{Ku}, Littelmann \cite{Li} and Mathieu \cite{M}.
In particular, it has been shown that in this case Schubert varieties are
normal and have rational singularities 
(see, for example, \cite{Ku}, theorems 2.16 and 2.23).

Now consider the affine flag variety over an algebraically closed
field of characteristic $p>0$. We want to show:

\begin{thm} {\bf (Faltings)} All Schubert varieties in $\cF$ are
normal, with rational singularities. Further, each Schubert variety
 $X\subset \cF$
is Frobenius split, and all Schubert subvarieties of $X$ are 
simultaneously compatibly split.
\end{thm}

For the notion of Frobenius splittings, which was introduced
by Mehta and Ramanathan, we refer to their articles
 \cite{MR} and \cite{Ram}. The facts that we will need can also be found
in \cite{G}.

The strategy of Faltings' proof is the following:

The Schubert varieties admit a resolution of singularities
 $ D(w) \lto C_w$
by the so-called Demazure varieties, which of course factors
over the normalization $\widetilde{C}_w$ of $C_w$.
One shows that the Demazure varieties are Frobenius split,
compatibly with their Demazure subvarieties, and that this
gives a Frobenius splitting of $\widetilde{C}_w$. 
Furthermore for $w' \le w$ we get a closed immersion
 $\widetilde{C}_w' \subseteq \widetilde{C}_w$. This allows us to define
the ind-scheme $\widetilde{\cC}$ of the normalizations 
 $\widetilde{C}_w$, and this ind-scheme can even be defined
over $\Z$. Now the key point is that the action of 
the  $SL_2(R) \subseteq G(R\xT)$ associated to the simple
affine roots can be lifted to an action on $\widetilde{\cC}$.
The subgroup $\cN^-$, which is the inverse image of the
unipotent radical of $B^-$ under the projection $L^-G \lto G$,
maps isomorphically onto an open subset of $\cF$.
Let $\cNN$ be the formal completion of $\cN^-$ at the origin.
We will see that $\cNN(R)$ is contained in the subgroup generated by the
 $SL_2(R)$'s associated to the simple affine roots. As a consequence
we get a section of the natural map between the completions
 of the local rings at the origin of $\widetilde{\cC}$ and $\cF$.
By exploiting the fact that the normality of Schubert varieties is
already known in characteristic 0, one shows that this
section is indeed an isomorphism; this implies the
normality result.

\subsection{Demazure varieties}

\subsubsection{Definitions}

In this section we define the Demazure varieties, which yield
a resolution of singularities of Schubert varieties.
We continue to work over an algebraically closed field $k$ and
we proceed exactly as in the case of the special linear group;
cf. \cite{M}, \cite{G}. Therefore the proofs are omitted.

Let $w \in W_a$ and let $ \tilde{w} = s_{i_1} \cdots s_{i_\ell}$ be
a reduced expression for $w$. (We sometimes
write $\tilde{w}$ instead of $w$ to
indicate that the following definitions really depend on the
choice of a reduced decomposition.)

\begin{Def} The variety
$$
D(\tilde{w}) := \cP_{i_1} \times^\cB \cdots \times^\cB \cP_{i_\ell}/\cB.
$$
is called the Demazure variety 
corresponding to $\tilde{w}$.
\end{Def}

If $\tilde{u} = s_{i_1} \cdots \widehat{s_{i_k}} \cdots s_{i_\ell}$ 
is reduced, we have a closed immersion
\begin{equation} \label{closedimm}
\sigma \colon D(\tilde{u}) \lto
 D(\tilde{w}). 
\end{equation}

If $\tilde{w} = \tilde{u} \tilde{v}$, i. e. 
 $\tilde{u} = s_{i_1} \cdots s_{i_k}$, 
 $\tilde{v} = s_{i_{k+1}} \cdots s_{i_\ell}$, for some $k$,
 we get a canonical projection
morphism $D(\tilde{w}) \lto D(\tilde{u})$,
which is a locally trivial fibre bundle with fibre $D(\tilde{v})$.
In particular, let $\tilde{u} = s_{i_1} \cdots s_{i_{\ell-1}}$,
 $\tilde{v} = s_{i_\ell}$. Then we get a $\PP^1$-fibration
\begin{equation} \label{p1fibr}
\pi \colon D(\tilde{w}) \lto D(\tilde{u}). \end{equation}
The closed immersion
 $D(\tilde{u}) \lto D(\tilde{w})$ defined above is a section
of this fibration.
Multiplication gives us a morphism 
 $\Psi_{\tilde{w}} : D(\tilde{w}) \lto X_w$.

\begin{stz} \label{demvar_compatible}
i) The Demazure variety $D(\tilde{w})$ is smooth and proper over $k$,
and has dimension $l(w)$.

ii) 
The morphism $\Psi_{\tilde{w}}:D(\tilde{w}) \lto X_w$ 
is proper and birational.
If $\tilde{u} = s_{i_1} \cdots \widehat{s_{i_k}} \cdots s_{i_\ell}$ 
is reduced,
these morphisms together with the closed immersion
(\ref{closedimm}) yield a commutative diagram
\begin{diagram}
D(\tilde{u}) & \rTo & X_u \\
\dInto & & \dInto \\
D(\tilde{w}) & \rTo & X_w
\end{diagram}
\end{stz}

In the Demazure variety $D(\tilde{w})$, we have $l(w)$ divisors
 $Z_1^{\tilde{w}} ,\dots, Z_{l(w)}^{\tilde{w}} $. 
These are defined inductively on the
length of $w$, as follows.

Write $\tilde{w} = \tilde{u} s_{i_\ell}$. We have the map
 $\pi: D(\tilde{w}) \lto D(\tilde{u})$, which is a $\PP^1$-fibration,
and we also have the section $\sigma: D(\tilde{u}) \lto D(\tilde{w})$
of $\pi$, as defined above.

We define
\begin{equation} \label{def_divisors}
\begin{array}{rcl}
Z_i^{\tilde{w}} & := & \pi^{-1}(Z_i^{\tilde{u}}), \quad i = 1,\dots, l(w)-1 
   \\[.2cm]
Z_{l(w)}^{\tilde{w}} & := & \sigma(D(\tilde{u})).
\end{array}
\end{equation}

We denote by $Z^{\tilde{w}}$ the sum of the divisors $Z_i^{\tilde{w}}$.

\begin{lem} \label{prop_zi} \label{codim1var_appear}
i) The subvarieties $Z^{\tilde{w}}_i$ are smooth of codimension 1 in
 $D(\tilde{w})$. 

ii) The scheme-theoretic intersection
 $P^{\tilde{w}}:= \bigcap_i Z_i^{\tilde{w}}$ is just a point.
 
iii)
If $\tilde{v} < \tilde{w}$, $l(v) = l(w)-1$,
then $D(\tilde{v})$ (considered as a closed subscheme of
 $D(\tilde{w})$ by the embedding defined above) is one of
the $Z_i^{\tilde{w}}$. \qed
\end{lem}

\subsubsection{The canonical bundle}

Let $w \in W_a$ and choose a reduced
decomposition $\tilde{w}$.
We want to describe the canonical bundle of the Demazure
variety $D(\tilde{w})$. 
As above, we identify the affine flag variety with the space
of special 
complete selfdual lattice chains. 
Again we denote by $(\lambda_i)_i$
the standard lattice chain.

The Schubert variety $X_w$ consists of certain lattice chains
 $(\sL_i)_i$. We can find $N > 0$, such that all lattices
occuring here lie between $t^{-N} k\xt^n$ and $t^N k\xt^n$.

Let $n_i = \dim_k \lambda_i/t^N k\xt^n$, $i= 0, \dots, r$.
We get maps 
$$ \varphi_i \colon X_w \lto \Grass(t^{-N}k\xt^n/t^N k\xt^n, n_i), \
 (\sL_i)_i \mapsto \sL_i/t^N k\xt^n. $$ 
These maps yield a closed
embedding
\begin{equation} \label{emb_in_grass}
\begin{diagram}
\varphi\colon X_w & \rInto & \prod_{i=0}^{r} 
       \Grass(t^{-N}k\xt^n/t^N k\xt^n, n_i).
\end{diagram}
\end{equation}

Now let $L_i$ be the very ample generator of the Picard group
of $\Grass(t^{-N}k\xt^n/t^N k\xt^n, n_i)$, and define 
$$ L_w := \varphi^\ast \bigotimes_{i=0}^{n-1} L_i. $$
This line bundle does not depend on $N$.

Let us remark that we could also obtain this line bundle 
as the pull-back of a certain line bundle on the affine flag 
variety; see \cite{F}, section 2.




Denote by $L_{\tilde{w}}$ the pull back of $L_w$ along
the morphism $\Psi_{\tilde{w}}: D(\tilde{w}) \lto X_w$.
As in \cite{G}, one proves

\begin{stz}\label{deg_along_fibres}

i) Let $\tilde{u} < \tilde{w}$, $l(\tilde{u}) = \tilde{w} - 1$.
 The pull back of $L_{\tilde{w}}$ along the embedding
 $\sigma: D(\tilde{u}) \lto D(\tilde{w})$ is $L_{\tilde{u}}$.

ii) The line bundle $L_{\tilde{w}}$ does not have a base point.

iii) Write $\tilde{w} = \tilde{u} s_{i_\ell}$.
 The degree of $L_{\tilde{w}}$ along the fibres of 
 $\pi: D(\tilde{w}) \lto D(\tilde{u})$ is $1$. 
\end{stz}

Now the following description of the canonical bundle of 
 $D(\tilde{w})$ is a purely formal consequence of the
above (cf. \cite{M}, 
\cite{G}).

\begin{stz} \label{canbundle}
The canonical bundle of $D(\tilde{w})$ is
 $$ \omega_{D(\tilde{w})} = \cO(-Z^{\tilde{w}}) \tens L_{\tilde{w}}^{-1}. $$
\end{stz}

\subsubsection{Demazure varieties are Frobenius split}

Assume that our algebraically closed field $k$ has 
characteristic $p > 0$.
Now that we have computed the canonical bundle, we can apply
the criterion of Mehta and Ramanathan for Frobenius split 
varieties, and get 

\begin{stz} \label{demvarsplit}
The Demazure variety $D(\tilde{w})$ admits a Frobenius
splitting which compatibly splits
all the divisors $Z_i^{\tilde{w}}$. 
\end{stz}

\subsection{The ind-scheme $\widetilde{\cC}$}

At the moment, we still work over an algebraically closed field of 
characteristic $p > 0$.

Denote by $\Cnorm_w$ the normalization
of $C_w$. The resolution $\pi_w \colon D(w) \lto C_w$ factors
as 
\begin{diagram}
 D(w) & \rTo & \Cnorm_w & \rTo^{\psi_w} & C_w. 
\end{diagram}
If $w' \le w$, the map $D(w') \lto D(w)$ induces a map
 $\Cnorm_{w'} \lto \Cnorm_w$, which is
independent of the chosen resolution of $w$.

By push forward, the Frobenius splitting of $D(w)$ gives 
us a splitting of $\Cnorm_w$; all $\Cnorm_{w'}$, $w'<w$ 
are simultaneously compatibly split.

The map $\pi_w$ has geometrically connected fibres as can be seen
by induction on the length of $w$ and by factoring
 $\pi_w$ as 
\begin{diagram}
D(w) = \cP_i \times^{\cB} D(w') & \rTo & \cP_i \times^{\cB} C_{w'}
 & \rTo & C_w,
\end{diagram}
where $w = s_i w'$; indeed, the first map has geometrically 
connected fibres by induction, and the fibres of 
the second map are points respectively projective lines.

This implies that the normalization map $\psi_w$ is a 
universal homeomorphism.

\begin{lem}
The natural maps $\widetilde{C}_{w'} \lto \widetilde{C}_{w}$,
     $w'\le w$, are closed immersions.
\end{lem}

{\em Proof.} Since $\psi_w$ is a universal homeomorphism,
the maps $\widetilde{C}_{w'} \lto \widetilde{C}_{w}$ 
are universally injective.

Denote by $C'$ the image of $\widetilde{C}_{w'}$ in 
 $\widetilde{C}_{w}$. Then the extension of the function
fields of  $\widetilde{C}_{w'}$ and $C'$ is purely
inseparable.

On the other hand, these varieties are compatibly
Frobenius split, i.e. we have a commutative diagram
of $\cO_{C'}$-modules
\begin{diagram}
\cO_{C'} & \rTo & \cO_{\widetilde{C}_{w'}} \\
\dTo^\varphi & & \dTo^\psi \\
\cO_{C'} & \rTo & \cO_{\widetilde{C}_{w'}} 
\end{diagram}
where $\varphi$ resp. $\psi$ are sections to the Frobenius
map $x \mapsto x^p$. From this we get an analogous diagram for the 
function fields, and this implies that they are indeed isomorphic.
\qed

\begin{lem} \label{rat_sing}
The higher direct 
 images $R^i \pi_{w,\ast} \cO_{D(w)}$,
 $i>0$, vanish. Furthermore $\widetilde{C}_w$ has rational 
singularities and in particular is Cohen-Macaulay.
\end{lem}

{\em Proof.}
To show the first part of the lemma,
write $w = s_i u$, $l(u) = l(w)-1$. Let us factor $\pi_w$ as
$$ D(w) = \cP_i \times^\cB D(u) 
  \lto \cP_i \times^\cB \widetilde{C}_u \lto \widetilde{C}_w. $$
By induction, we may assume that the higher direct images of
 $\cO_{D(w)}$ under the first morphism vanish. To deal
with the second part, observe that its fibre over $x \in \widetilde{C}_w$
is either a point or a $\PP^1$. Thus we can apply the following lemma 
(see \cite{MS} for a proof) and are done.

\begin{lem} {\bf (Mehta - Srinivas)}
Let $f : X \lto Y$ be a projective birational morphism of algebraic
varieties over $k$, such that
$X$ is Frobenius split, and for all $y \in Y$ we have $H^i(X_y, \cO_{X_y})
=0$ for all $i > 0$. Then $R^i f_\ast \cO_X = 0$ for all $i > 0$.\qed
\end{lem}

Finally, to show that $\pi_w$ is a rational resolution, 
it remains to show that the higher direct images of the canonical
bundle vanish as well. This follows from the Grauert-Riemenschneider
theorem for Frobenius split varieties, see \cite{MvK}.
\qed

Now we define the ind-scheme $\widetilde{\cC}$ as
the inductive limit of the normalizations $\Cnorm_w$, $w \in W_{a}$.

We even can do this over $\Z$; to do so, we will use the following
proposition which is a nice generalized version of the well-known
cohomology and base change theorem.

\begin{stz} {\bf (Faltings)} \label{coh_bc} 
Let $S$ be a noetherian scheme, let $X$ and $Y$ be noetherian
 $S$-schemes,
and let $f \colon X \lto Y$ be a proper morphism of
$S$-schemes.
Furthermore, assume that $F$ is a coherent sheaf on $X$
which is flat over $S$. If $R^1 f_{s,\ast} F_s =0$ for all
 $ s\in S$, then $f_\ast F$ is flat over $S$ and is compatible
with base change $S' \lto S$.
\end{stz} 

{\em Proof.} 
It is clear that $f_\ast F$ is flat over $S$.
To show the compatibility with base change,
we more or less follow section III.12 in
Hartshorne's book \cite{H}.

Without loss of generality, we may assume
that $S$ and $Y$ are affine, say $S=\Spec A$,  $Y=\Spec B$.

Now consider the (covariant) functor
$$ T^i(M) = H^i(X, F \otimes_A M) $$
on the category of $A$-modules.
The functor $T^i$ is additive and exact in the middle. The family
 $(T^i)_{i \ge 0}$ is a $\delta$-functor.

\begin{lem} (see {\bf \cite{H} III 12.5})
For any $M$, we have a natural map $\varphi: T^i(A) \otimes M \lto T^i(M)$.
Furthermore, the following conditions are equivalent:

i) $T^i$ is right exact

ii) $\varphi$ is an isomorphism for all $M$

iii) $\varphi$ is surjective for all $M$.
\end{lem}

{\em Proof.} The proof in loc. cit. carries over literally. \qed

Let us first show that $T^1(M)=0$ for every $A$-module $M$. 
By a flat base change, we may assume that $A$ is 
local, with maximal ideal $\fm$. We
denote the closed point of its spectrum by $s$. 
Furthermore, it is enough to show the claim for finitely
generated $A$-modules, since $T^i$ commutes with direct
limits.

First consider the case of a module of finite
length.
By our assumption $T^1(k(s))=0$, i. e. $T^1$
vanishes on the $A$-module of length one. 
Since $T^1$ is exact in the middle,
this implies that $T^1(M)=0$ for all $A$-modules $M$ of
finite length. 

Now let $M$ be an arbitrary $A$-module of finite type.
Assume that $T^1(M)$ does not vanish, and let $y \in \Spec B$
be a point in the support of $T^1(M)$. By further
localizing $A$ if necessary, we may assume that $y$ maps
to the closed point $s \in \Spec A$.

We know already that $T^1(M/\fm^n M)=0$ for all $n$.
So by the theorem on formal functions, applied to $f:X\lto Y$,
we get 
$$ 0= \lim_{\longleftarrow} T^1(M/\fm^n M) = H^1(X, F\otimes M)^\wedge, $$
where $H^1(X, F\otimes M)^\wedge$ denotes the completion of the $B$-module
 $H^1(X, F \otimes M)$ with respect to the ideal $\fm B$.
But since $\fm B \subseteq \fp_y$, this implies that the stalk $T^1(M)_y$
vanishes which is a contradiction.

So we have proved that $T^1$ vanishes, and thus
that $T^0$ is right exact. By the lemma above
this yields that $T^0(M) \cong T^0(A) \otimes M$ for all $A$-modules
 $M$, which is what we needed to prove. \qed

Now we can define the ind-scheme of normalizations of
Schubert varieties over the integers.
The Demazure varieties $D(w)$
are obviously defined over $\Z$. Now we 
let $\Cnorm_w = \Spec \pi_\ast(\cO_{D(w)})$, where $\pi$
is the map $D(w) \lto \cF$ induced by multiplication.
We get flat $\Z$-schemes that are independent of the chosen
reduced decomposition of $w$ (since $\Spec \pi_\ast(\cO_{D(w)})$
is the normalization of the scheme-theoretic image of
$D(w)$ in $\cF$).
This construction is
compatible with base change by lemma \ref{rat_sing}
and proposition \ref{coh_bc}.
Furthermore, $\Cnorm_{w'}$ is a closed subscheme of $\Cnorm_w$
for $w' \le w$, and we obtain an ind-scheme
$$ \widetilde{\cC} = \lim_{\lto} \Cnorm_w $$
over $\Z$, which maps to $\cF$.

\begin{lem} The action of the $SL_2$'s 
associated to the simple affine roots lifts to $\widetilde{\cC}$.
\end{lem}

{\em Proof.} Fix $i \in \{0, \dots, r \} $.
The subgroup $SL_2$ associated to a simple affine root
 $\alpha_i$ acts on $C_w$ if $l(s_iw) = l(w)-1$. 
If $l(s_iw) = l(w)-1$, we can even find a reduced
decomposition of $w$ that starts with $s_i$. Then 
 $SL_2$ acts on the Demazure variety corresponding to
this decomposition (since it acts on the
parahoric subgroup $\cP_i$).
Thus we also get an $SL_2$-action  on $\Cnorm_w$.

Since there
is a cofinal system of elements $w$ with $l(s_iw) = l(w)-1$ 
in $W_{a}$, the action lifts to $\widetilde{\cC}$.
\qed

\subsection{The open cell}

We continue to work over the integers. Let us for a moment
consider the special linear group $SL_n$. Again, we denote
by $B$ the standard Borel subgroup, by $B^-$ the opposite
Borel subgroup, and by $\cB$ the Iwahori subgroup of $L\,SL_n$
corresponding to $B$.

Denote by $\cN^- = \cN^-_{SL_n}$ the inverse image 
of the unipotent radical of $B^-$
under the projection $L^- SL_n = SL_n(k[t^{-1}]) \lto SL_n$.

\begin{stz}
The canonical map $\cN^- \lto L\,SL_n \lto L\,SL_n/\cB$ is 
an open immersion.
\end{stz}

Before we prove the proposition, let us first look at the
situation in the affine Grassmannian, which is similar.
Faltings (\cite{F}, lemma 2) shows that the natural map 
 $L^{--}SL_n \lto L\, SL_n / L^+ SL_n$ is an open immersion,
the image of which consists of those lattices $\sL \subset R\xT^n$,
such that
$$ \sL \oplus t^{-1} R[t^{-1}]^n = R\xT^n. $$
If $\sL$ satisfies this condition, 
we get the corresponding matrix in $L^{--}SL_n$ in the following way:
Denote by $e_1, \dots, e_n$ the standard basis of $R\xT^n$, and write
$$ e_i = u_i + m_i, \quad u_i \in \sL,\ m_i \in t^{-1} R[t^{-1}]^n. $$
Then it can be shown that the $u_i$ are a basis for $\sL$,
and the matrix $(u_1, \dots, u_n) = 1- (m_1, \dots, m_n)$ lies in
 $L^{--}SL_n$.

Now let us generalize this observation to the affine flag variety,
and thus prove the proposition.

{\em Proof.} Let $R$ be a $\Z$-algebra. Denote
the standard lattice chain by $(\lambda_i)_i$. Furthermore, write 
 $\sM_i = t^{-2}R[t^{-1}]^i \oplus t^{-1}R[t^{-1}]^{n-i}$.
We will show that the image of $\cN^-(R)$ is the set of lattice chains
 $(\sL_i)_i$ such that 
\begin{equation} \label{cond}
 \sL_i \oplus \sM_i = R\xT, \qquad\text{ for all } i. 
\end{equation}
This is an open condition; more precisely, we see that $\cN^-$
is identified with the intersection of the inverse images 
of the open cells of the affine Grassmannians.

It is clear that the condition (\ref{cond}) is satisfied by
lattice chains of the form $h((\lambda_i)_i)$, $h\in \cN^-(R)$. Now 
take a lattice chain $(\sL_i)_i \in LSL_n/\cB$ that 
satisfies (\ref{cond}). 

For each $i \in \{0,\dots,n-1\}$ write
\begin{eqnarray*}
 t^{-1} e_j & = & u^i_j + m^i_j, \quad u^i_j \in \sL_i,
                                 \ m^i_j \in \sM_i,\ j=1, \dots, i \\
        e_j & = & u^i_j + m^i_j, \quad u^i_j \in \sL_i,
                                 \ m^i_j \in \sM_i,\ j=i+1, \dots, n
\end{eqnarray*}
where again $e_1, \dots, e_n$ denotes the standard basis of 
 $R\xT^n$.
By what we have said above, the $u^i_j$, $j=1,\dots, n$,
are a basis of $\sL_i$.

Claim: The matrix $h:=(tu^1_1, tu^2_2, \dots, tu^{n-1}_{n-1}, u^0_n)$
lies in $\cN^-(R)$, and $h(\lambda_i)_i = (\sL_i)_i$.

By the definition of the $u^i_j$, it is clear that $h\in \cN^-(R)$.
So it remains to show that $h(\lambda_i) = \sL_i$, for all $i$.
This is done by an easy computation, and to keep the notations
simple, we will only treat the case $i=0$.
We have, for $1 \le j < n$,
\begin{eqnarray*}
e_j & = & t(u^j_j + m^j_j) \\
    & = & t u^j_j + \sum_{k=j+1}^n a_k e_k
           + m^0, \quad a_k \in R,\ m^0 \in \sM_0,  \\
    & = & t u^j_j + \sum_{k=j+1}^n a_k u^0_k + m'_0, \quad m'_0 \in \sM_0,
\end{eqnarray*}
so we have 
$$ u^0_j = t u^j_j + \sum_{k=j+1}^n a_k u^0_k,$$
and this proves, by descending induction on $j$, that indeed
 $tu^1_1, \dots, tu^{n-1}_{n-1}, u^0_n$ is a basis of $\sL_0$.
\qed

Now we switch back to the symplectic group $G = Sp_{2r}$. 
Inside $L^-G$, we have the subgroup $\cN^-_{Sp_{2r}}$, which is
by definition the inverse image of the unipotent radical
 $N^-$ of the opposite Borel subgroup $B^-$, under
the projection $L^-G \lto G$.
In fact, $\cN^-_{Sp_{2r}}$ is just the intersection (inside $L\,SL_n$)
of $\cN^-_{SL_{2r}}$ with $LG$.
From now on, we will simply write $\cN^-$ instead of $\cN^-_{Sp_{2r}}$.
As before, we denote the affine flag variety $L\,Sp_{2r}/\cB$ by $\cF$.
The proposition gives us the following corollary.

\begin{kor} The natural map $\cN^- \lto \cF$ is an open immersion.
\qed
\end{kor}

We denote by $\cNN$ the formal completion of $\cN^-$
at the identity. In other words,
\begin{eqnarray*}
 \cNN(R) = \{ g \in G(R[t^{-1}]); && g \equiv 1 \text{ modulo 
     a nilpotent ideal in $R$,} \\
 && \text{ and the constant term of $g$ 
     lies in $N^-$}  \}. 
\end{eqnarray*}
It is clear (think of the elements of $\cNN$ as matrices)
that  $\cNN$ is ind-represented by a power series ring 
 $\cR := \Z[\hspace{-.12em}[ x_i, i\in \N ]\hspace{-.12em} ]$ 
in infinitely many variables,
i. e. 
$$ \cNN = \lim_{\lto} \Spec \cR/I, $$
where the limit runs over those ideals $I$ which contain
almost all $x_i$, and contain some power of each $x_i$.
These ideals $I$ will be called open.
It is clear that we can just as well take the limit over
those $\Spec \cR/I$, such that additionally $\cR/I$ is
 $\Z$-torsion free.

\begin{lem}
Let $g \in \cNN(R)$. Then the image of $g$ in $G(R\xT)$ lies in the
subgroup $K$ generated by the subgroups $SL_2(R)$ associated to the
simple affine roots.
\end{lem}

{\em Proof.}
Of course, the subgroups $U_\alpha$ associated to the simple affine roots $\alpha$
are contained in $K$.
On the other hand, we get all elements of the affine Weyl group, 
since the simple reflections are contained in the $SL_2(R)$'s.
Thus all $U_\alpha$ lie in $K$.

The entries on the diagonal of $g$ have the form $1 + t^{-1} a$,
 $a \in {\rm nil}(R)$, hence are units in $R[t^{-1}]$.
Thus we can write $g$ as a product of 
elements of $L^{--}T$ and of certain $U_\alpha$'s.

So it only remains to show that $L^{--}T$ is contained
in $K$.
But obviously $L^{--}T$ is generated by the subgroups $\G_m(R[t^{-1}])
 \cong \{ \diag(1^i, a, 1^{2r-2i-2}, a^{-1}, 1^i) \}$
associated to the (finite) simple roots.
To show that these $\G_m(R[t^{-1}])$ are contained in $K$,
it suffices to treat the case of $SL_2$. In this case the assertion 
follows from the formula ($a\in R[t^{-1}]^\times$)
$$ 
\left(\begin{array}{cc}
a & 0 \\
0 & a^{-1}
\end{array} \right)
=
\left(\begin{array}{cc}
1 & a \\
0 & 1
\end{array} \right)
\left(\begin{array}{cc}
1 & 0 \\
-a^{-1} & 1
\end{array} \right)
\left(\begin{array}{cc}
1 & a \\
0 & 1
\end{array} \right)
\left(\begin{array}{cc}
0 & -1 \\
1 & 0
\end{array} \right).
$$
\qed

\subsection{Proof of the normality theorem} 

Denote by 
 $\cA$ the projective limit of the local rings 
of the formal completions of $\Cnorm_w$ at the origin.
Since $\widetilde{\cC}$ maps to $\cF$, we have a natural map
 $\cR \lto \cA$. We will first define a section of this
map, and then show that it is an isomorphism.

Let $I \subseteq \cR$ be an open ideal such that $R:= \cR/I$
is torsion free over $\Z$, and let $g \in \cNN(R)$.
By the lemma in the previous section, $g$ lies in the
subgroup of $LG$ generated by $SL_2(R)$'s, so we can
lift the corresponding point $z \in \cF$ to a point  
 $\tilde{z} \in \widetilde{\cC}$, say $\tilde{z} \in \Cnorm_{w}$,
which again lies in the formal completion at the origin.

On the level of rings, this means the following:  The map $\cR \lto R$
corresponding to $g$ lifts to a continous map $\cA \lto R$.
This lift is unique because $R$ is $\Z$-torsion free
and $\Cnorm_{w,\Q} \lto C_{w,\Q}$ is an isomorphism.
By taking the limit, we get a section $\cA \lto \cR$
of the natural map $\cR \lto \cA$.

We want to show that these maps actually are isomorphisms,
i. e. that the composition $\cA \lto \cR \lto \cA$
is the identity map as well.

Let $w \in W_{a}$, and denote by $A$ the ring of the formal
completion of $\Cnorm_w$ at the origin. Choose a lift of
the origin to a $\Z$-valued point of the Demazure variety
 $D(w)$, and denote the ring of the formal completion
at this point by $B$. The ring $B$ is simply a ring of 
formal power series over $\Z$. By $B_n$ we will 
denote the quotient of $B$ modulo some power of its
defining ideal.

We have an inclusion  $A \subseteq B$,
and  it is enough to show that the two endomorphisms
of $\cA$, namely the composition $\cA \lto \cR \lto \cA$ and
the identity map,
 become equal after composing them with $\cA \lto A \lto B \lto B_n$.
Now these two compositions correspond to $B_n$-valued
points of $\Cnorm_w$ which become equal in $\cF$. Since 
 $B_n$ is torsion-free and $\Cnorm_{w,\Q} \cong C_{w,\Q}$,
these points are actually equal.

So in fact $\cR \cong \cA$, and this implies that 
the map $\Cnorm_w \lto \cF$ is a closed immersion
near the origin, because it induces a surjection
on the local rings of the formal completions.

Now consider the situation over an algebraically closed field $k$
again. What we have shown implies that $C_w$ is normal
at the origin. On the other hand, the locus where $C_w$ is not
normal is closed and $\cB$-stable, hence, if it was not empty, it 
would have to contain the origin, which is the unique closed 
 $\cB$-orbit.

\section{Doubly symplectic standard tableaux}
\label{deConcini}

Let $r \ge 1$ be an integer, $n=2r$, let $k$ be a field and
$$ R= k[c_{\mu \nu}, \mu, \nu = 1, \dots 2r] / (CJ_{2r}C^t, C^tJ_{2r}C).
$$
Here $C=(c_{\mu\nu})_{\mu,\nu}$ and $J_{2r}$ again denotes the matrix 
 $\left( \begin{array}{cc} & {\mathbf J} \\ -{\mathbf J} & 
\end{array} \right)$, where 
 ${\mathbf J}$ is the $r \times r$-matrix with
1's on the second diagonal and 0's elsewhere. 

We will see later that this singularity appears in the 
'simplest' singular local models, and it will be important for
us to know that $R$ is reduced. To prove reducedness, we make
use of a theorem of de Concini \cite{dC} which we explain now.

The theorem will give a $k$-basis of $R$ in terms of 
doubly symplectic standard tableaux; to give the definition
of these standard tableaux,
we need several notations. 

First of all, if $I=\{i_1, \dots, i_k \}$ and $J=\{ j_1, \dots, j_\ell \}$
are subsets of $\{ 1, \dots, r\}$, we say $ I \le J$
if $ k \ge \ell $ and $i_\mu \le j_\mu$, $\mu = 1, \dots, \ell$.

Now if $i_1, \dots, i_s, j_1, \dots, j_s \in \{ 1, \dots, r\}$,
we denote by $(i_s, \dots, i_1 | j_1, \dots, j_s)$ the minor
of $C$ consisting of the rows $i_1, \dots, i_s$ and columns
 $j_1, \dots, j_s$.

Since we are in a sense dealing with the symplectic group, it will be useful
to consider the rows (resp. columns) $1, \dots, r$ and $r+1, \dots, 2r$
separately. Thus we make the following definition:

Let $s \le k \le r$, 
$I = \{ i_1, \dots, i_s \}, J = \{ j_1, \dots, j_{k-s} \}
 \subseteq \{ 1, \dots, r \}$. Let $\Gamma = I \cap J$,
and write $I = \tilde{I} \stackrel{\cdot}{\cup} \Gamma$,
 $J = \tilde{J} \stackrel{\cdot}{\cup} \Gamma$. Finally,
let $\tilde{I} = \{ \tilde{i}_1 < \dots < \tilde{i}_{s-\lambda} \}$,
 $\tilde{J} = \{ \tilde{j}_1 < \dots < \tilde{j}_{k-s-\lambda} \}$,
 $\Gamma = \{ \gamma_1 < \dots < \gamma_\lambda \}$.
Then we define 
\begin{eqnarray*}
 (J, I | h_1, \dots, h_k ) & = &
   (n-\tilde{j}_1+1, \dots, n-\tilde{j}_{k-s-\lambda}+1, 
    \tilde{i}_{s-\lambda}, \dots, \tilde{i}_1, \\
&&    n-\gamma_\lambda+1, \gamma_\lambda, \dots, n-\gamma_1+1, \gamma_1 
      | h_1, \dots, h_k ).  
\end{eqnarray*}

Next we have to define the notion of admissible minor.
So let $I$, $J$ and $\Gamma = I\cap J$ be as above.

\begin{Def}
The minor $P = (J, I | h_1, \dots h_k)$ is
called admissible, if there exists a subset
 $T \subseteq \{1, \dots, r \} - ( I \cup J)$, such that 
 $ |T| = |\Gamma|$ and $T \ge \Gamma$.

In this case we will write 
\begin{eqnarray*}
 P & = & \left(\left. \begin{array}{c} {J', I} \\ {J, I'} \end{array}
 \right| h_1, \dots, h_k \right) \\
 & = & \left(\left. \begin{array}{c} n-j_1'+1, \dots, n-j_{k-s}'+1,
          i_s, \dots, i_1 \\ n-j_1+1, \dots, n-j_{k-s}+1,
          i_s', \dots, i_1' \end{array}
 \right| h_1, \dots, h_k \right),
\end{eqnarray*}
where $I = \tilde{I} \cup \Lambda = \{ i_1' < \dots < i_s' \}$, 
 $J' = \tilde{J} \cup \Lambda = \{ j_1' < \dots < j_{k-s}' \} $,
and where $\Lambda \subseteq \{1, \dots, r \} - ( I \cup J)$
is the smallest subset with the above properties.
\end{Def}

{\em Example.} The minor $(r, \dots, 1 | 1, \dots, r) = 
(\emptyset, \{1, \dots, r\} | 1, \dots, r)$
is obviously admissible. We have $\Gamma = \Lambda = \emptyset$,
and
$$ (\emptyset, \{1, \dots, r\} | 1, \dots, r)
 =  \left(\left. \begin{array}{c} {r, \dots, 1 } \\
                    { r, \dots, 1 } \end{array}
     \right| 1, \dots, r \right).
$$

The minor $P$ will be called doubly admissible, if the admissibility
condition holds not only for the rows, but also for the columns.
In this case, we can write
\begin{eqnarray*} P &  = & 
        \left(\left. \begin{array}{c} J', I \\ J, I' \end{array}
        \right| \begin{array}{c} H, K' \\ H', K \end{array} \right) \\
   & = &    \left(\left. \begin{array}{c} t_k, \dots, t_1 \\ 
                                          v_k, \dots, v_1 \end{array}
 \right|  \begin{array}{c} s_1, \dots, s_k \\ 
                           u_1, \dots, u_k \end{array}  \right)            .
\end{eqnarray*}

We have some kind of partial order on the set of doubly admissible
minors, that is defined in the following way.
Let 
$$ P = \left(\left. \begin{array}{c} t_k, \dots, t_1 \\ 
                           v_k, \dots, v_1 \end{array}
 \right|  \begin{array}{c} s_1, \dots, s_k \\ 
             u_1, \dots, u_k \end{array}  \right), 
\quad 
P' = \left(\left. \begin{array}{c} t'_{k'}, \dots, t'_1 \\ 
                                   v'_{k'}, \dots, v'_1 \end{array}
 \right|  \begin{array}{c} s'_1, \dots, s'_{k'} \\ 
                           u'_1, \dots, u'_{k'} \end{array}  \right)   $$
be doubly admissible minors. Then we say that $P \le P'$
if 
$$
\{ v_1, \dots, v_k \} \le \{ t'_1, \dots t'_{k'}\}
\text{ and }
\{ u_1, \dots, u_k \} \le \{ s'_1, \dots s'_{k'}\}.
$$
Note that, in general, it is not true that $P \le P$.

\begin{Def}
A tuple $(P_1, \dots, P_\ell)$ of doubly admissible minors,
such that $P_i \le P_{i+1}$ for all $i$, is
called a doubly symplectic standard tableaux.
We denote the set of all doubly symplectic standard tableaux
by $\cT$.
\end{Def}

We have a map $\varphi \colon \cT \lto R$, which maps an element
 $(P_1, \dots, P_\ell)$ of $\cT$ to the product $P_1 \cdot \dots P_\ell$
of the minors $P_i$.

With these notations, de Concini (\cite{dC}, theorem 6.1) 
proved the following

\begin{thm} {\bf (de Concini)}
 The elements $\varphi(P)$, $P \in \cT$, form a $k$-basis of $R$.
\end{thm}

We will in fact only need the following corollary to de Concini's theorem.

\begin{kor} \label{f_NNT}
The element $f = (r, \dots, 1 | 1, \dots, r)$
is not a zero divisor in $R$.
\end{kor}

{\em Proof.} Obviously, $f$ is a doubly admissible minor; see the
example above. Furthermore, the set $\{1, \dots, r\}$ is the smallest 
element among all subsets of $\{1, \dots, r\}$ with respect to the partial
order defined above. Thus, if $T$ is an arbitrary doubly symplectic
standard tableau, $fT$ is again a doubly symplectic standard tableau.
In view of the theorem this clearly implies that $f$ is not a zero 
divisor.

{\bf Remark.}
What we really need to know, and what we will derive from this corollary
in the course of the proof of the flatness theorem, is that
 $R$ is reduced; see proposition \ref{Rred}. 
It might well be that this result is - at least
implicitly - already contained in de Concini's paper.
For example, it would follow if one knew that $R$
equipped with the basis of doubly symplectic standard tableaux
becomes an algebra
with straightening law. I do not know if this is true.

\section{Proof of the flatness theorem}
\label{proof}

We embed the special fibre $\Mlocs$ of the local model 
into the affine flag variety
$\cF = LG/\cB$. 
More precisely, we identify $\cF$ with the space
of $r$-special selfdual complete lattice chains, as explained
in proposition \ref{F_latticechains} and the remark following it.
Let $R$ be a $k$-algebra. As before, we denote by $(\lambda_i)_i$
the standard lattice chain:
$$ \lambda_i = t^{-1} R\xt^i \oplus R\xt^{n-i}. $$
Then we map an element
 $(\cF_i)_i \in \Mlocs(R)$ to
the lattice chain $(\sL_i)_i$, where $\sL_i$ is the
inverse image of $\cF_i$, considered as a subspace of 
 $\Lambda_{i,R} \cong R^n$,
under the projection
 $\lambda_i \lto \lambda_i/t \lambda_i \cong R^n$. 
 
It is easily checked that in this way we get an element of
 $\cF(R)$, and that the resulting map is a closed immersion.
Via this embedding, we can identify $\Mlocs$
with the set
\begin{equation} \label{descrMloc}
 \Mlocs \cong \{ (\sL_i)_i \in \cF; \ t \lambda_i \subseteq \sL_i \subseteq \lambda_i
 \text{ for all } i \}. 
\end{equation}

We can do the same thing for the parahoric local model associated
to a subset $I \subseteq \{ 1, \dots, n \}$ (such that $2r-i \in I$
for each $i\in I$, $i\ne 0$). Denote by $\cP^I \subseteq LG$ 
the stabilizer of the partial lattice chain $(\lambda_i)_{i\in I}$.
We get an embedding $\Mlocs_I \subset LG/\cP^I$.
Denote by 
 $\Mloct_I$ the inverse image of $\Mlocs_I \subset Sp_{2r}(k\xT)/P^I$
under the projection $\cF \lto Sp_{2r}(k\xT)/P^I$.
It can be identified with the following set of lattice chains:
$$ \Mloct_I = \{ (\sL_i)_i \in \cF; \ t \lambda_i \subseteq 
\sL_i \subseteq \lambda_i \text{ for all } i \in I \}. $$

Then we have (scheme-theoretic intersection)
$$ \Mlocs = \Mloct_0 \cap \Mloct_r \cap 
           \bigcap_{i=1}^{r-1} \Mloct_{i, 2r-i}. $$

All the $\Mloct_I$ and also $\Mlocs$ are, set-theoretically,
unions of Schubert varieties in $\cF$, since they
are invariant under the action of $\cB$.




\subsection{The simplest local models}
\label{simplest}

The local models $\Mloc_0$ and $\Mloc_{r}$ are isomorphic to
the Grassmannian of totally isotropic subspaces of $\Lambda_0$
resp. $\Lambda_r$. So they are smooth over $O$.

The case of $\Mloc_{i,2r-i}$, $i \in \{ 1, \dots, r-1 \}$, is more
complicated. First we determine the equations of the local model in this
case. For simplicity we assume that $r \le 2i$ (the case $r > 2i$ is 
completely analogous).

In an open neighborhood of the worst singularity, we can represent 
the subspaces $\cF_i, \cF_{2r-i}$ with 
$(\cF_i, \cF_{2r-i}) \in \Mloc_{i,2r-i}$ by matrices of the form

$$
\cF_i \hat{=}
\left(
\begin{array}{ccc}
a_{r-i+1,1} & \cdots & a_{r-i+1,r} \\
\vdots      &        &  \vdots     \\
a_{r1}      & \cdots & a_{rr}      \\
1           &        &             \\
            & \ddots &             \\
            &        &    1        \\
 a_{11}     & \cdots &  a_{1r}     \\
\vdots      &        &  \vdots     \\
a_{r-i,1}   & \cdots & a_{r-i,r}   
\end{array}
\right), \qquad
\cF_{2r-i} \hat{=}
\left(
\begin{array}{cccccc}
        &        &        &      1     &        &             \\
        &        &        &            & \ddots &             \\
        &        &        &            &        &    1        \\
 b_{11} & \cdots & b_{1i} & b_{1,i+1}  & \cdots &   b_{1r}     \\
\vdots  &        & \vdots &            & \vdots &    \vdots    \\
b_{r1}  & \cdots & b_{ri} & b_{r,i+1}  & \cdots &  b_{rr}     \\
  1     &        &        &            &        &             \\
        & \ddots &        &            &        &             \\
        &        &    1   &            &        &                
\end{array}
\right).
$$

The duality condition implies that 
$$b_{\mu \nu} = \varepsilon_{\mu \nu} a_{r-\nu+1, r-\mu+1},$$
where 
$$\varepsilon_{\mu \nu} = 
 \left\{ \begin{array}{rl} 1 & \mu,\nu \le i \text{ or } 
 \mu, \nu \ge i+1, \\ -1 & \text{otherwise.} \end{array}\right.
$$

To write down in terms of the $a_{\mu \nu}$ what it means that these
subspaces are mapped into one another, we divide the matrices
 $ A = (a_{\mu \nu})_{\mu, \nu} $ and $B = (b_{\mu \nu})_{\mu, \nu}$
into blocks:
$$ 
\begin{array}{rcr}
A = &
\left(
\begin{array}{cc}
A_4 & A_3 \\
A_1 & A_2
\end{array}
\right) & \begin{array}{c} \mbox{\tiny $2(r-i)$} \\ \mbox{\tiny $2i-r$} 
\end{array} \\
& \begin{array}{cc} \mbox{\tiny $2(r-i)$} & \mbox{{\tiny $2i -r$}} 
\end{array} &  
\end{array}, \qquad
\begin{array}{rcr}
B = &
\left(
\begin{array}{cc}
B_2 & B_3 \\
B_1 & B_4
\end{array}
\right) & \begin{array}{c} \mbox{\tiny $ 2i-r$} \\ \mbox{\tiny $2(r-i)$} 
\end{array} \\
& \begin{array}{cc} \mbox{\tiny $2i-r$} & \mbox{{\tiny $ 2(r-i)$}}
\end{array} &
\end{array}.
$$
Now we get the following equations as condition that $\cF_i$ 
is mapped into $\cF_{2r-i}$:
$$ B_4 A_4 = \pi, \quad B_2 = A_2 - B_3 A_3, \quad A_1 = B_3 A_4, \quad
   B_1 = - B_4 A_3.
$$
The condition that $\cF_{2r-i}$ is mapped into $\cF_i$ is expressed 
by the following equations:
$$ A_4 B_4 = \pi, \quad \pi B_2 = A_1 B_1 + \pi A_2, \quad
   \pi B_3 = A_1 B_4, \quad \pi A_3 = -A_4 B_1. $$
Altogether, this gives us the following:
$$ B_4 A_4 = A_4 B_4 = \pi, \quad B_2 = A_2 - B_3 A_3, \quad A_1 = B_3
A_4. $$
All the other equations follow from these four equations. 

It is easy to see that we can choose all the coefficients of $A_3$ and
those of $A_2$ that lie on or above the secondary diagonal. Then the
other coefficients of $A_2$ are determined, as well as the coefficients
of $A_1$. Of course, $B_1$, $B_2$ and $B_3$ are then determined, too.

Furthermore
note that $B_4$ is just the adjoint matrix of $A_4$ with respect
to the standard symplectic form.

So our open subset is isomorphic to
\begin{eqnarray*}
&& \Spec O[c_{\mu \nu}, \mu, \nu = 1, \dots 2(r-i)] / 
    (C C^{ad} = C^{ad} C = \pi) \\
&&\times 
    \A_O^{2(i-r)(2i-r)+ \frac{1}{2}(2i-r)(2i-r+1)}.
\end{eqnarray*}

Thus to show that $\Mlocs_{i,2r-i}$ is reduced, we may obviously
restrict ourselves to the case where $r = 2i$.
Since $C^{ad} = -J_{2i} C^t J_{2i}$, where $J_{2i}$ is the matrix
describing the standard symplectic form as above,
 in this case an open neighborhood 
of the worst singularity in the special fibre is isomorphic to
$\Spec R$, where
$$  R= k[c_{\mu \nu}, \mu, \nu = 1, \dots r] / 
    (CJ_{2i}C^t, C^tJ_{2i}C),
$$
and we have to show 

\begin{stz} \label{Rred} The ring $R$ is reduced.
\end{stz}

{\em Proof.}
We denote by $f\in R$ the $i \times i$-minor of $C$ consisting of 
the rows $1, \dots, i$ and columns $1, \dots, i$, and by $R_f$
the localization of $R$ with respect to $f$.

\begin{lem} The scheme $\Spec R_f$ is smooth; more precisely,
it is an open subscheme of $\A_k^{r(r+1)/2}$.
\end{lem}

{\em Proof.} Consider the open subset
 $U \subset \Mlocs_{i,2r-i}$, where we can represent the subspaces as
$$ \cF_i = 
\left(
\begin{array}{ccc}
a_{11} & \cdots & a_{1r} \\
\vdots & & \vdots \\
a_{r1} & \cdots & a_{rr} \\
1 & & \\
 & \ddots & \\
 & & 1
\end{array}
\right), \qquad
\cF_{2r-i} = 
\left(
\begin{array}{ccc}
b_{11} & \cdots & b_{1r} \\
\vdots & & \vdots \\
b_{r1} & \cdots & b_{rr} \\
1 & & \\
 & \ddots & \\
 & & 1
\end{array}
\right).
$$
The duality condition implies that $b_{\mu \nu} = a_{r-\nu+1,, r-\mu+1}$.

It is easy to see by a direct computation
that this open subset is isomorphic
to $\A_k^{r(r+1)/2}$. 
We could also derive this by using the corresponding fact for $GL_n$;
cf. section
\ref{liftability}.

Since $\Spec R_f$ obviously is an open subscheme of $U$, 
the lemma follows.
\qed

Furthermore, corollary \ref{f_NNT} shows that $f$ is not
a zero divisor. 

Since $\Spec R_f$ is irreducible and $f$ is not a zero divisor,
we see that $\Mloc_{i, 2r-i}$ is irreducible.
This can also be proved by analysing the affine Weyl group
(cf. the article \cite{KR} of Kottwitz and Rapoport, sections 4 and 10).

Now we can conclude that $R$ is reduced: obviously this is true generically,
so all we have to show is that $R$ has no embedded associated 
prime ideals. This is clear for the localization $R_f$,
and on the other hand
all primes in the complement of $\Spec R_f$ contain
the non-zero divisor $f$, and thus cannot be associated prime ideals.
\qed

\subsection{Intersections of Schubert varieties}
\label{inters}

From Faltings' theorem we get 

\begin{kor} Let $Y_1 , \dots Y_k \subset \cF$ be unions
of Schubert varieties. Then the scheme-theoretic intersection
 $\bigcap_i Y_i$ is reduced.
\end{kor}

{\em Proof.} Embed all the $Y_i$ in a sufficiently big Schubert
variety $X$, and use the fact that $X$ is Frobenius split, compatibly
with all its Schubert subvarieties. \qed

The corollary shows that $\Mlocs$ is 
a union of Schubert varieties in the affine flag variety 
even scheme-theoretically.
In particular, it is reduced and its irreducible components
are normal with rational singularities.

The final ingredient for the flatness theorem is that we can lift the
generic points of the special fibre to the generic fibre.
This point is dealt with in the next section.

\subsection{Lifting the generic points of the special fibre}
\label{liftability}

We want to show that the generic points of the special
fibre of a local model $\Mloc_I$ can be lifted to the generic fibre.
To keep the notation a little bit simpler, we deal only with 
 $\Mloc$. See \cite{G} for a detailed account on this
question in the case of $GL_n$.

As we have seen, $\Mlocs$ is the union of certain Schubert
varieties, that correspond to certain elements of the affine
Weyl group, or in other words to certain alcoves (in the standard
apartment of the Bruhat-Tits building) for $Sp_{2r}$. 

The question which Schubert varieties occur has a natural answer
in terms of the notions of admissible and permissible alcoves
introduced by Kottwitz and Rapoport \cite{KR},
which we will use freely in the following. 

It is clear that the set of alcoves contributing to $\Mlocs$
is precisely the set of minuscule or $\mu$-permissible alcoves, where 
 $\mu$ is the minuscule dominant coweight $(1^r, 0^r)$.
As Kottwitz and Rapoport have shown, this set coincides with the
set of $\mu$-admissible alcoves, which consists of all alcoves
that are smaller than some conjugate of $\mu$ under the finite Weyl group.

Thus the 'extreme' elements, which correspond to the irreducible
components of $\Mlocs$, are the conjugates of $\mu$ under $W$.

Let $x = (x_0, \dots, x_{2r-1})$ be an extreme minuscule alcove for
 $GSp_{2r}$ (cf. \cite{KR}). Then $x_i = x_0 + (1^i, 0^{2r-i})$,
and $\{ x_0(i), x_0(2r-i+1) \} = \{ 0, 1 \}$, $i = 1, \dots, r$.

Let $I = \{ i_1 < \cdots < i_r \} = \{ i; x_0(i) = 0 \}$,
 $J = \{ j_1 < \cdots < j_r \} = \{ j; x_0(j) = 1 \}$.

We denote by $U_x$ the open subset of $\Mloc$, where each subspace
 $\cF_i$ can be represented by a matrix with the 
unit matrix in rows $j_1, \dots, j_r$. Note that for each
 $\cF_i$ the unit matrix is in the same place.
The special fibre $\overline{U}_x$ of $U_x$ contains the Schubert cell
corresponding to $x$ as an open subscheme.

\begin{stz} We have $U_x \cong \A_O^{r(r+1)/2}$.
\end{stz}

{\em Proof.}
We can consider $U_x$ as a closed subscheme of the corresponding
extreme stratum $U_x(GL_{2r})$ in the local model for $GL_{2r}$,
and we know that $U_x(GL_{2r}) \cong \A^{r^2}$ (cf. \cite{G}).

Recall how this isomorphism was constructed. We denote
the coefficients (in the rows $i_1, \dots, i_r$)
in the matrix describing $\cF_i$ by $a^i_{\lambda \mu}$.

Now fix $\lambda$ and $\mu$. Then for $i=0, \dots, 2r-1$
we get three types of equations:

{\em First case.} $\pi a^i_{\lambda \mu} = a^{i+1}_{\lambda \mu}.$

{\em Second case.} $a^i_{\lambda \mu} = \pi a^{i+1}_{\lambda \mu}.$

{\em Third case.} $a^i_{\lambda \mu} = a^{i+1}_{\lambda \mu}.$

The first and second case occur precisely once, namely when
 $i = i_\lambda - 1$ resp. $i = j_\mu$.
So we see that we can choose $a^{i_\lambda-1}_{\lambda \mu}$
arbitrarily, and that all other $a^i_{\lambda \mu}$ are then uniquely
 determined; more precisely we have 
$$ a^i_{\lambda \mu} = 
\left\{ 
\begin{array}{rl} 
a^{i_\lambda-1}_{\lambda \mu} & \quad i \in [ j_\mu, i_\lambda ) \\
\pi a^{i_\lambda-1}_{\lambda \mu} & \quad\text{otherwise}
\end{array}
\right. ,
$$
where we consider $[ j_\mu, i_\lambda )$ as an interval in
 $\Z/2r\Z$.

Now we have to investigate the effect of the duality condition
that describes $U_x$ inside $U_x(GL_{2r})$.
It says that for all $i, \lambda, \mu$, we must have
 $a^i_{\lambda \mu} = a^{2r-i}_{r-\mu+1, r-\lambda+1}$.
We want to get our isomorphism $U_x \cong \A^{r(r+1)/2}$
by choosing just one of these two, and the
following lemma shows that we can do so.

\begin{lem} For all $\lambda, \mu$, we have
$$ i\in [j_\mu, i_\lambda )\quad \Longleftrightarrow \quad
   2r-i \in [j_{r-\lambda+1}, i_{r-\mu+1}).
$$
\end{lem}

{\em Proof.} Since we started with an extreme minuscule
alcove {\em for the symplectic group},
we have
$$ i_\lambda = 2r - j_{r-\lambda+1} +1, \quad
   j_\mu = 2r- i_{r-\mu+1} +1,
$$
and the lemma follows immediately.
\qed

Since an open subset of the irreducible component of $\Mlocs$
corresponding
to $x$ is contained in $U_x$ as an open subscheme, the proposition
implies that the generic point of this irreducible
component can be lifted to the generic fibre.
This concludes the proof of the flatness theorem.

\vskip1cm
\hspace*{1cm}Ulrich G\"ortz\\
\hspace*{1cm}Mathematisches Institut\\
\hspace*{1cm}der Universit\"at zu K\"oln\\
\hspace*{1cm}Weyertal 86--90

\hspace*{1cm}{\small DE--}50931 K\"oln (Germany)

\hspace*{1cm}ugoertz@mi.uni-koeln.de

\vskip.5cm
\hspace*{1cm}{\em until july 2001:}\\
\hspace*{1cm}School of Mathematics\\
\hspace*{1cm}Institute for Advanced Study\\
\hspace*{1cm}Einstein Drive

\hspace*{1cm}Princeton, NJ 08540 (USA)

\hspace*{1cm}ugoertz@ias.edu
\end{document}